\documentclass[12pt]{article}
\usepackage{amssymb,amsmath}
\begin{document}
\font\SY=msam10
\def\N{{\mathbb N}}
\def\Z{{\mathbb Z}}
\def\R{{\mathbb R}}
\def\C{{\mathbb C}}
\def\T{{\mathbb T}}
\def\zp{\Z_+}
\def\dist{{\rm dist}\,}
\def\Sim{{\rm Sim}\,}
\def\square{\hbox{\SY\char"03}}
\def\epsilon{\varepsilon}
\def\phi{\varphi}
\def\kappa{\varkappa}
\def\theorem#1{\smallskip\noindent
{\scshape Theorem} {\bf #1}{\bf .}\hskip 8pt\sl}
\def\defin#1{\smallskip\noindent
{\scshape Definition} {\bf #1}{\bf .}\hskip 6pt}
\def\prop#1{\smallskip\noindent
{\scshape Proposition} {\bf #1}{\bf .}\hskip 8pt\sl}
\def\lemma#1{\smallskip\noindent
{\scshape Lemma} {\bf #1}{\bf .}\hskip 6pt\sl}
\def\cor#1{\smallskip\noindent
{\scshape Corollary} {\bf #1}{\bf .}\hskip 6pt\sl}
\def\quest#1{\smallskip\noindent
{\scshape Question} {\bf #1}{\bf .}\hskip 6pt\sl}
\def\epr{\smallskip\rm}
\def\rem#1{\smallskip\noindent
{\scshape Remark}  {\bf #1}.\hskip 6pt}
\def\wz{\thinspace}
\def\proof{P\wz r\wz o\wz o\wz f.\hskip 6pt}
\def\leq{\leqslant}
\def\geq{\geqslant}
\def\ssub#1#2{#1_{{}_{{\scriptstyle #2}}}}
\def\H{{\cal H}}
\def\B{{\cal B}}
\def\slim{\mathop{\hbox{$\overline{\hbox{\rm lim}}$}}}
\def\ilim{\mathop{\hbox{$\underline{\hbox{\rm lim}}$}}}

\title{Compact operators without extended eigenvalues}

\author{S.~Shkarin}

\date{}

\maketitle

\smallskip

\leftline{King's College London, Department of Mathematics}
\leftline{Strand, London WC2R 2LS, UK} \leftline{\bf e-mail: \tt
stanislav.shkarin@kcl.ac.uk}

\rm\normalsize

\begin{abstract}
A complex number $\lambda$ is called an extended eigenvalue of a
bounded linear operator $T$ on a Banach space $\B$ if there exists a
non-zero bounded linear operator $X$ acting on $\B$ such that
$XT=\lambda TX$. We show that there are compact quasinilpotent
operators on a separable Hilbert space, for which the set of
extended eigenvalues is the one-point set $\{1\}$.
\end{abstract}

\section{Introduction}

All vector spaces in this article are over the field $\C$ of complex
numbers. For a Banach space $\B$, $L(\B)$ stands for the algebra of
bounded linear operators on $\B$. A complex number $\lambda$ is
called an {\it extended eigenvalue} of $T\in L(\B)$ if there exists
non-zero $X\in L(\B)$ such that $XT=\lambda TX$. We denote the set
of extended eigenvalues of $T$ by the symbol $\Sigma(T)$. Extended
eigenvalues and their corresponding extended eigenoperators $X$ as
well as their applications to classification of invariant subspaces
of linear operators were studied in \cite{blp,bp,lam}. Obviously
$1\in\Sigma(T)$ for any operator $T$. Indeed, one can take $X$ being
the identity operator.

Let $V$ be the Volterra operator acting on $L_2[0,1]$:
$$
V:L_2[0,1]\to L_2[0,1],\quad Vf(x)=\int\limits_0^1 f(t)\,dt.
$$
In \cite{blp} the set of extended eigenvalues of $V$ has been
computed.

\theorem{BLP}The set $\Sigma(V)$ coincides with the set $\R^+$ of
positive real numbers. \epr

In \cite{lam} it is also shown that $\Sigma(\gamma I+V)=\{1\}$ for
each non-zero $\gamma\in\C$. Note that $\gamma I+V$ is neither
compact nor quasinilpotent. The following questions were raised in
\cite{blp}. In what follows $\H$ stands for a separable infinite
dimensional Hilbert space.

\quest1 Does there exist a compact operator $T\in L(\H)$ for which
$\Sigma(T)=\{1\}?$ \epr

\quest2 Does there exist a quasinilpotent operator $T\in L(\H)$ for
which $\Sigma(T)=\{1\}?$ \epr

We answer both questions affirmatively:

\theorem{1.1}There exists a compact quasinilpotent operator $T\in
L(\H)$ for which $\Sigma(T)=\{1\}$. \epr

It is worth noting that $\Sigma(T)=\C$ if $T$ is a nilpotent
operator acting on a Banach space $\B$. Indeed, if $T=0$, then
$XT=\lambda TX=0$ for each $X\in L(\B)$ and each $\lambda\in\C$. If
$T\neq 0$, there exists a positive integer $n$ such that $T^n\neq 0$
and $T^{n+1}=0$. Then $XT=\lambda TX=0$ is satisfied with $X=T^n$
for any $\lambda\in\C$.

The proof of Theorem~1.1 is non-constructive. The existence of a
required operator $T$ is proved by using the Baire category
argument. The proof is based on a theorem of Apostol, characterizing
closures of similarity orbits of compact quasinilpotent operators
and precise knowledge of the sets of extended eigenvalues for two
specific operators, one of which is the Volterra operator and the
other one is a bilateral weighted shift.

We compute extended eigenvalues of bilateral weighted shifts and
prove other auxiliary results in Section~2. Theorem~1.1 is proved in
Section~3. In Section~4 of concluding remarks we discuss previous
results and raise few problems.

\section{Auxiliary results}

From now on $\Z$ stands for the set of integers and $\zp$ is the set
of non-negative integers.

\prop{2.1}Let $1\leq p<\infty$, $w=\{w_n\}_{n\in\Z}$ be a bounded
sequence of non-zero complex numbers and $T$ be the bilateral
weighted shift with the weight sequence $w$ acting on the space
$\ell_p(\Z)$, that is $Te_n=w_ne_{n-1}$ for any $n\in\Z$, where
$\{e_n\}_{n\in\zp}$ is the canonical basis of $\ell_p(\Z)$. Let also
\begin{equation}
\beta(k,n)=\prod_{j=k}^n w_j\ \ \text{for}\ \ k,n\in\Z,\ k\leq n\ \
\text{and \ $\beta(n+1,n)=1$ \ for $n\in\Z$}. \label{beta}
\end{equation}
Then $0\notin\Sigma(T)$ and a non-zero complex number $\lambda$
belongs to $\Sigma(T)$ if and only if there exists $k\in\zp$ such
that the sequence $\{\lambda^{-n}\beta(n-k+1,n)\}_{n\in\Z}$ is
bounded. \epr

\proof Clearly $T$ is injective and has dense range. Therefore zero
is not an extended eigenvalue of $T$. First, assume that there
exists $k\in\zp$ for which the sequence
$\{\lambda^{-n}\beta(n-k+1,n)\}_{n\in\Z}$ is bounded. Then there
exists a unique bounded linear operator $X$ on $\ell_p(\Z)$ such
that $Xe_n=\lambda^{-n}\beta(n-k+1,n)e_{n-k}$ for any $n\in\Z$.
Clearly $X\neq 0$. It is straightforward to verify that $XT=\lambda
TX$. Hence $\lambda\in\Sigma(T)$.

Assume now that the sequence
$\{\lambda^{-n}\beta(n-k+1,n)\}_{n\in\Z}$ is unbounded for any
$k\in\zp$.  Let $X\in L(\ell_p(Z))$ be such that $XT=\lambda TX$. We
have to prove that $X=0$. Let $\{x_{n,j}\}_{n,j\in\Z}$ be the matrix
of $X$, that is a complex matrix such that
$$
Xe_n=\sum_{j\in\Z}x_{n,j} e_j\ \ \text{for any $n\in\zp$}.
$$
From the equation $XT=\lambda TX$ it immediately follows that
$w_nx_{n-1,j}=\lambda w_{j+1}x_{n,j+1}$ for any $j,n\in\Z$.
Therefore
\begin{equation}
x_{n,n+k}=\left\{\begin{array}{ll}
\lambda^{-n}\frac{\beta(1,k)}{\beta(n+1,n+k)}x_{0,k}&\text{if
$k>0$}; \vrule width0pt height0pt depth12pt
\\
\lambda^{-n}\frac{\beta(k+n+1,n)}{\beta(k+1,0)}x_{0,k}&\text{if
$k\leq0$}\end{array} \right. \label{xnk} \qquad \text{for any
$n,k\in\Z$.}
\end{equation}
Boundedness of $X$ implies boundedness of the set
$\{x_{n,j}:n,j\in\Z\}$. If there exists $k\leq 0$ for which
$x_{0,k}\neq 0$, then according to (\ref{xnk}), boundedness of
$\{x_{n,n+k}\}_{n\in\Z}$ implies boundedness of
$\{\lambda^{-n}\beta(k+n+1,n)\}_{n\in\Z}$, which is unbounded by
assumption. If there exists $k>0$ for which $x_{0,k}\neq 0$, then
according to (\ref{xnk}) boundedness  of $\{x_{n,n+k}\}_{n\in\Z}$
implies boundedness of $\{\lambda^{-n}/\beta(n+1,n+k)\}_{n\in\Z}$.
Since $\{w_n\}_{n\in\Z}$ is bounded, the sequence
$\{\beta(n+1,n+k)\}_{n\in\Z}$ is also bounded, and it follows that
$\{\lambda^{-n}=\lambda^{-n}\beta(n+1,n)\}_{n\in\Z}$ is bounded as a
product of two bounded sequences, which also contradicts the
assumption. Thus, $x_{0,k}=0$ for each $k\in\Z$. Formula (\ref{xnk})
now implies that $x_{n,j}=0$ for any $n,j\in\Z$ and therefore $X=0$.
\square

\rem1 Let for $k\in\zp$,
$$
c_k=\slim_{n\to+\infty}\bigl|\beta(n+1,n+k)\bigr|^{1/n}\ \ \text{and}\ \
d_k=\ilim_{n\to+\infty}\bigl|\beta(1-n,k-n)\bigr|^{-1/n}.
$$
Then $c_{k+1}\leq c_k\leq 1\leq d_k\leq d_{k+1}$ for each $k\in\zp$.
Denote $c=\lim\limits_{k\to\infty}c_k$ and
$d=\lim\limits_{k\to\infty}d_k$. Then $0\leq c\leq 1\leq
d\leq\infty$. From Proposition~2.1 it follows that $\Sigma(T)$ is an
annulus of one of the following shapes $\{z\in\C:c\leq |z|\leq d\}$,
$\{z\in\C:c< |z|\leq d\}$, $\{z\in\C:c\leq |z|<d\}$ or
$\{z\in\C:c<|z|< d\}$. Moreover, $\Sigma(T)$ always contains the
unit circle.

\cor{2.2}Let $\H$ be a separable infinite dimensional Hilbert space.
Then there exists an injective quasinilpotent compact operator $T$
on $\H$ with dense range such that $\Sigma(T)$ coincides with the
unit circle $\T=\{z\in\C:|z|=1\}$. \epr

\proof Since all separable infinite dimensional Hilbert spaces are
isomorphic, we can assume that $\H=\ell_2(\Z)$. Let $T:\H\to\H$ be
the bilateral weighted shift with the weight sequence
$w_n=(1+|n|)^{-1}$. Then we have $\beta(n-k+1,n)\sim (1+|n|)^{-k}$
as $n\to\infty$ for each $k\in\zp$, where the numbers $\beta(a,b)$
are defined by (\ref{beta}). Hence for any non-zero complex number
$\lambda$ and any $k\in\zp$, the sequence
$\{\lambda^{-n}\beta(n-k+1,n)\}_{n\in\Z}$ is bounded if and only if
$|\lambda|=1$. According to Proposition~2.1, $\Sigma(T)=\T$. Clearly
$T$ is compact, injective and has dense range. Since for any
$k\in\zp$ and any $n\in\Z$, $T^ke_n=\beta(n-k+1,n)e_{n-k}$, we have
\begin{equation}
\|T^k\|=\max_{n\in\Z}\beta(n-k+1,n). \label{e1}
\end{equation}
Now, if $\alpha_n=\beta(n-k+1,n)$, then according to (\ref{beta}),
$$
\frac{\alpha_{n}}{\alpha_{n-1}}=
\frac{w_n}{w_{n-k}}=\frac{1+|n-k|}{1+|n|}.
$$
Therefore $\frac{\alpha_{n}}{\alpha_{n-1}}\leq 1$ if $n\geq k/2$ and
$\frac{\alpha_{n}}{\alpha_{n-1}}\geq 1$ if $n\leq k/2$. Thus,
$\alpha_n$ is maximal  for $n=[k/2]$, where $[k/2]$ is the integer
part of $k/2$. Hence
$$
\max_{n\in\Z}\beta(n-k+1,n)=\max_{n\in\Z}\alpha_n=
\alpha_{[k/2]}=\left\{\begin{array}{ll}
(m!)^{-2},\ \text{where $m=(k+1)/2$}&\text{if $k$ is odd};\\
(m!(m+1)!)^{-1},\ \text{where $m=k/2$}&\text{if $k$ is even.}
\end{array}\right.
$$
The above display together with (\ref{e1}) implies
$\lim\limits_{k\to\infty}\|T^k\|^{1/k}=0$. Thus, from the spectral
radius formula \cite{rud} it follows that $T$ is quasinilpotent.
\square

The following lemma is a direct consequence of Theorem~BLP.

\lemma{2.3}Let $\H$ be a separable infinite dimensional Hilbert
space. Then there exists an injective quasinilpotent compact
operator $T$ on $\H$ with dense range such that $\Sigma(T)$
coincides with the open half-line $\R^+=\{x\in\R:x>0\}$. \epr

For a Banach space $\B$, symbol $\Omega(\B)$ stands for the set of
compact quasinilpotent operators $T\in L(\B)$ endowed with the
operator norm metric $d(T,S)=\|T-S\|$. It is easy to show that
$(\Omega(\B),d)$ is a complete metric space, or equivalently, that
$\Omega(\B)$ is closed in the Banach space $K(\B)$ of compact
operators on $\B$. Indeed, if $T\in K(\B)$ non-quasinilpotent, then
$T$ has a non-zero normal eigenvalue. Hence all operators
sufficiently close to $T$ with respect to the metric $d$ also have a
non-zero normal eigenvalue and therefore can not be quasinilpotent.
It is worth noting that we actually need completeness of
$\Omega(\B)$ only in the case when $\B=\H$ is a Hilbert space. In
this case it can be found, for instance, in \cite{herre}.

For a bounded linear operator $T$ acting on a Banach space $\B$ we
denote
\begin{equation}
\Sim(T)=\{GTG^{-1}:G\in L(\B)\ \ \text{is invertible.}\}
\label{sim}
\end{equation}
The set $\Sim(T)$ is usually called the {\it similarity orbit} of
$T$. The key tool for the proof of Theorem~1.1 is the following
theorem of Apostol, which is the main result in \cite{apo}. It worth
noting that Apostol's theorem answers a question raised by Herrero
in \cite{her}, where a weaker statement was proved.

\theorem{A}Let $\H$ be a separable infinite dimensional Hilbert
space and $T\in\Omega(\H)$ be non-nilpotent. Then $\Sim(T)$ is dense
in $\Omega(\H)$. \epr

Recall that a subset $A$ of a topological space $Y$ is called an
$F_\sigma$-set if $A$ is a union of countably many closed subsets of
$Y$. The complement of an $F_\sigma$-set is called a $G_\delta$-set.
Equivalently a $G_\delta$-set in $Y$ is a countable intersection of
open subsets of $Y$.

\lemma{2.4}Let $\B$ be a separable reflexive Banach space and
$A\subset\C$ be an $F_\sigma$-set. Then
$$
M(A)=\{T\in L(\B):\Sigma(T)\cap A\neq\varnothing\}
$$
is an $F_\sigma$-set in $L(\B)$ endowed with the topology defined by
the operator norm metric $d(T,S)=\|T-S\|$. \epr

\proof Obviously any $F_\sigma$-set in $\C$ is a union of countably
many compact sets. Pick compact sets $K_n\subset\C$, $n\in\zp$ such
that $A=\bigcup\limits_{k=0}^\infty K_n$. Since $\B$ is separable,
there exists a sequence $\{x_m\}_{m\in\zp}$ of elements of $\B$ such
that the set $\{x_m:m\in\zp\}$ is dense in $\B$. Reflexivity of $\B$
implies that $\B^*$ is also separable and we can choose a sequence
$\{y_k\}_{k\in\zp}$ of elements of $\B^*$ such that
$\{y_k:k\in\zp\}$ is dense in $\B^*$. For $n,m,k,j\in\zp$ we denote
$$
N(n,m,k,j)=\biggl\{T\in L(B):\begin{array}{l} \text{there exist
$X\in L(\B)$ and $\lambda\in K_n$ such that}\\
\text{$XT=\lambda TX$, $\|X\|\leq j$ and $\langle Xx_m,y_k\rangle
=1$}\end{array}\biggr\}.
$$
First, we shall show that
\begin{equation}
M(A)=\bigcup_{n,m,k,j\in\zp} N(n,m,k,j). \label{union1}
\end{equation}
The inclusion $N(n,m,k,j)\subseteq M(A)$ follows immediately from
the definitions of $N(n,m,k,j)$ and $M(A)$. Let $T\in M(A)$. Then
there exists $\lambda\in A$ and non-zero $Y\in L(\B)$ such that
$YT=\lambda TY$. Since $A=\bigcup\limits_{k=0}^\infty K_n$, there
exists $n\in\zp$ such that $\lambda\in K_n$. Since $Y\neq 0$ and
$\{x_m:m\in\zp\}$ is dense in $\B$, there exists $m\in\zp$ for which
$Yx_m\neq 0$. Since $\{y_k:k\in\zp\}$ is dense in $\B^*$, there
exists $k\in\zp$ such that $\langle Yx_m,y_k\rangle\neq 0$. Let
$X=(\langle Yx_m,y_k\rangle)^{-1}Y$. Then $XT=\lambda TX$ and
$\langle Xx_m,y_k\rangle=1$. Finally taking $j\in\zp$ such that
$j\geq\|X\|$, we ensure that $T\in N(n,m,k,j)$, which proves
(\ref{union1}).

In view of (\ref{union1}), it suffices to verify that the sets
$N(n,m,k,j)$ are closed with respect to the operator norm topology.
Let $n,m,k,j\in\zp$, $T\in L(\B)$ and $\{T_l\}_{l\in\zp}$ be a
sequence of elements of $N(n,m,k,j)$ such that $\|T_l-T\|\to 0$ as
$l\to\infty$. We have to verify that $T\in N(n,m,k,j)$. Since
$T_l\in N(n,m,k,j)$, there exist $\lambda_l\in K_n$ and $X_l\in
L(\B)$ such that $T_lX_l=\lambda_l X_lT_l$, $\|X_l\|\leq j$ and
$\langle X_lx_m,y_k\rangle =1$ for each $l\in\zp$. Since $\B$ is
separable and reflexive for any $c>0$ and sequence
$\{Z_l\}_{l\in\zp}$ of elements of $L(\B)$ such that $\|Z_l\|\leq c$
for each $l\in \zp$, there exists $Z\in L(\B)$ with $\|Z\|\leq c$
and a strictly increasing sequence  $\{l_q\}_{q\in\zp}$ of positive
integers such that $Z_{l_q}$ converges to $Z$ as $q\to\infty$ with
respect to the weak operator topology \cite{rud}, that is $Z_{l_q}x$
converges weakly to $Zx$ for each $x\in\B$. Using this fact along
with compactness of $K_n$, we can, passing to a subsequence, if
necessary, assume that $\lambda_l\to\lambda\in K_n$ and $X_l\to X\in
L(\B)$ as $l\to\infty$ with respect to the weak operator topology,
where $\|X\|\leq j$. Since $\langle X_lx_m,y_k\rangle =1$ for each
$l\in\zp$, we see that
$\langle Xx_m,y_k\rangle =1$. Let $x\in\B$. Then
\begin{equation}
(T_lX_l-\lambda_l X_lT_l)x=(\lambda-\lambda_l)X_lT_lx+
(T-T_l)X_lx+\lambda X_l(T-T_l)x+TX_lx-\lambda X_lTx.
 \label{equ}
\end{equation}
Since the sequence $T_l$ is norm-convergent, there exists $c>0$ such
that $\|T_l\|\leq c$ for each $l\in\zp$. Then
\begin{alignat}{2}
\|(\lambda-\lambda_l)X_lT_lx\|&\leq jc\|x\||\lambda-\lambda_l|\to
0&&\text{as $l\to\infty$}, \notag
\\
\|(T-T_l)X_lx\|&\leq j\|x\|\|T-T_l\|\to 0&&\text{as $l\to\infty$},
\notag
\\
\|\lambda X_l(T-T_l)x\|&\leq j|\lambda|\|x\|\|T-T_l\|\to 0\qquad
&&\text{as $l\to\infty$}. \notag
\end{alignat}
Since $X_lx$ converge weakly to $Xx$ and any bounded linear operator
on a Banach space is weak-to-weak continuous, we see that $TX_lx$
converge weakly to $TXx$. Finally $X_lTx$ converge weakly to $XTx$.
Thus, from (\ref{equ}) and the last display it follows that
$(T_lX_l-\lambda_l X_lT_l)x$ converge weakly to $(TX-\lambda XT)x$.
On the other hand $(T_lX_l-\lambda_l X_lT_l)x=0$ for any $l\in\zp$
and therefore $(TX-\lambda XT)x=0$ for each $x\in\B$. Thus,
$TX-\lambda XT=0$, which proves the inclusion $T\in N(n,m,k,j)$.
\square

\section{Proof of Theorem~1.1}

Let $A_1=\C\setminus \R^+$ and $A_2=\C\setminus\T$. Clearly $A_1$
and $A_2$ are $F_\sigma$-sets in $\C$. Let $\H$ be a separable
infinite dimensional Hilbert space and
$$
M_j=\{T\in\Omega(\H):\Sigma(T)\cap A_j\neq\varnothing\},\quad j=1,2.
$$
According to Lemma~2.4 $M_j$ are $F_\sigma$-sets in $\Omega(\H)$. By
Lemma~2.3 and Corollary~2.2, there exist non-nilpotent operators
$T_1,T_2\in\Omega(\H)$ such that $\Sigma(T_1)=\R_+$ and
$\Sigma(T_2)=\T$. Since the set $\Sigma(T)$ is a similarity
invariant, we have that $\Sigma(S)=\R_+$ for each $S\in\Sim(T_1)$
and $\Sigma(S)=\T$ for each $S\in\Sim(T_2)$. Hence
$$
\Sim(T_1)\cap M_1=\Sim(T_2)\cap M_2=\varnothing.
$$
By Theorem~A, both $\Sim(T_1)$ and $\Sim(T_2)$ are dense in
$\Omega(\H)$. According to the last display, $\Omega(\H)\setminus
M_1$ and $\Omega(\H)\setminus M_2$ are both dense $G_\delta$ sets in
$\Omega(\H)$. Since $\Omega(\H)$ is a complete metric space, the
Baire theorem implies that
$$
(\Omega(\H)\setminus M_1)\cap (\Omega(\H)\setminus M_2)=
\Omega(\H)\setminus (M_1\cup M_2)
$$
is a dense $G_\delta$-set in $\Omega(\H)$ and therefore $M_1\cup
M_2$ is a Baire first category set, that is a countable union of
nowhere dense sets. From the definitions of $M_1$ and $M_2$ and
the equality $\R^+\cap\T=\{1\}$ it
immediately follows that
$$
M_1\cup M_2=\{T\in\Omega(\H):\Sigma(T)\neq\{1\}\}.
$$
Thus, we have proven the following theorem.

\theorem{3.1}Let $\H$ be a separable Hilbert space. Then
$\{T\in\Omega(\H):\Sigma(T)\neq\{1\}\}$ is a Baire first category
set in the complete metric space $\Omega(\H)$. \epr

Theorem~1.1 is an immediate consequence of Theorem~3.1 and the Baire
theorem.

\section{Concluding remarks}

Theorem~BLP shows that the set $\Sigma(T)$ for a bounded linear
operator on a separable Hilbert space can be non-closed. The
following proposition shows that the sets $\Sigma(T)$ can not be too
bad. The proof goes similarly to the proof of Lemma~2.4. It is
included for sake of completeness.

\prop{4.1}Let $\B$ be a separable reflexive Banach space and $T\in
L(\B)$. Then $\Sigma(T)$ is an $F_\sigma$-set. \epr

\proof As in the proof of Lemma~2.4, we can choose sequences
$\{x_m\}_{m\in\zp}$ of elements of $\B$ and $\{y_k\}_{k\in\zp}$ of
elements of $\B^*$ such that $\{x_m:m\in\zp\}$ is dense in $\B$ and
$\{y_k:k\in\zp\}$ is dense in $\B^*$. For $m,k,j\in\zp$ we denote
$$
A(m,k,j)=\biggl\{\lambda\in\C:\begin{array}{l} \text{there exist
$X\in L(\B)$ such that $XT=\lambda TX$,}\\
\text{$\|X\|\leq j$ and $\langle
Xx_m,y_k\rangle=1$}\end{array}\biggr\}.
$$

First, we shall show that
\begin{equation}
\Sigma(T)=\bigcup_{m,k,j\in\zp} A(m,k,j). \label{union2}
\end{equation}
The inclusion $A(m,k,j)\subseteq \Sigma(T)$ is obviously valid for
any $m,k,j\in\zp$. Let $\lambda\in \Sigma(T)$. Then there exists
non-zero $Y\in L(\B)$ such that $YT=\lambda TY$. Since $Y\neq 0$ and
$\{x_m:m\in\zp\}$ is dense in $\B$, there exists $m\in\zp$ for which
$Yx_m\neq 0$. Since $\{y_k:k\in\zp\}$ is dense in $\B^*$, there
exists $k\in\zp$ such that $\langle Yx_m,y_k\rangle\neq 0$. Let
$X=(\langle Yx_m,y_k\rangle)^{-1}Y$. Then $XT=\lambda TX$ and
$\langle Xx_m,y_k\rangle=1$. Finally taking $j\in\zp$ such that
$j\geq\|X\|$, we ensure that $\lambda\in A(m,k,j)$, which proves
(\ref{union2}).

In view of (\ref{union2}), it suffices to verify that the sets
$A(m,k,j)$ are closed in $\C$. Let $m,k,j\in\zp$, $\lambda\in\C$ and
$\{\lambda_l\}_{l\in\zp}$ be a sequence of elements of $A(m,k,j)$,
converging to $\lambda$. Since $\lambda_l\in A(m,k,j)$, there exist
$X_l\in L(\B)$ such that $TX_l=\lambda_l X_lT$, $\|X_l\|\leq j$ and
$\langle X_lx_m,y_k\rangle =1$ for each $l\in\zp$. Since $\B$ is
separable and reflexive, we, passing to a subsequence, if necessary,
assume that $X_l\to X\in L(\B)$ as $l\to\infty$ with respect to the
weak operator topology, where $\|X\|\leq j$. Since $\langle
X_lx_m,y_k\rangle =1$ for each $l\in\zp$, we see that $\langle
Xx_m,y_k\rangle =1$. Let $x\in\B$. Then
\begin{equation}
(TX_l-\lambda_l X_lT)x=(\lambda-\lambda_l)X_lTx+ TX_lx-\lambda
X_lTx.
 \label{equ1}
\end{equation}
Clearly
$$
\|(\lambda-\lambda_l)X_lTx\|\leq j|\lambda-\lambda_l|\|T\|\|x\|\to
0\ \ \text{as $l\to\infty$}.
$$
Since $X_lx$ converge weakly to $Xx$ and any bounded linear operator
on a Banach space is weak-to-weak continuous, we see that $TX_lx$
converge weakly to $TXx$. Finally $X_lTx$ converge weakly to $XTx$.
Thus, from (\ref{equ1}) and the last display it follows that
$(TX_l-\lambda_l X_lT)x$ converge weakly to $(TX-\lambda XT)x$. On
the other hand $(TX_l-\lambda_l X_lT)x=0$ for any $l\in\zp$ and
therefore $(TX-\lambda XT)x=0$ for each $x\in\B$. Hence $TX-\lambda
XT=0$, which proves the inclusion $\lambda\in A(m,k,j)$. \square

Proposition~4.1 is a first step towards the solution of the
following problem.

\quest{4.2}Which subsets of $\C$ have the form $\Sigma(T)$ for $T\in
L(\H)$? Which subsets of $\C$ have the form $\Sigma(T)$ for compact
$T\in L(\H)$? Which subsets of $\C$ have the form $\Sigma(T)$ for
compact quasinilpotent subsets of $\C$? What about operators on
Banach spaces? \epr

The key point of the above proof of Theorem~1.1 is application of
Theorem~A, which does not work in the Banach space setting. This leads
naturally to the following question.

\quest{4.3}For which Banach spaces $\B$ there exists a compact
quasinilpotent operator $T\in L(\B)$ such that $\Sigma(T)=\{1\}?$
\epr

\bigskip

{\bf Acknowledgements.} \ Partially supported by British Engineering
and Physical Research Council Grant GR/T25552/01. The author is
grateful to Professor Alan Lambert for interest and helpful
comments.

\end{document}